\documentclass[a4paper]{article}

\usepackage{amsmath,amsfonts,amssymb,amsthm}
\usepackage[hidelinks,pdfusetitle]{hyperref}
\usepackage[capitalize]{cleveref}
\usepackage[a4paper, width=135mm, height=215mm]{geometry}

\numberwithin{equation}{section}
\newtheorem{theorem}{Theorem}
\newtheorem{proposition}{Proposition}[section]
\newtheorem{open}[proposition]{Open problem}
\newtheorem{definition}[proposition]{Definition}
\newtheorem{corollary}[proposition]{Corollary}
\newtheorem{lemma}[proposition]{Lemma}
\newtheorem{remark}[proposition]{Remark}

\newcommand{\supp}{\operatorname{supp}}
\newcommand{\N}{\mathbb{N}}
\newcommand{\R}{\mathbb{R}}
\newcommand{\dd}{\,\mathrm{d}}

\title{A family of interpolation inequalities involving products of low-order derivatives}

\author{Frédéric Marbach\texorpdfstring{\thanks{DMA, École normale supérieure, Université PSL, CNRS, 75005 Paris, France}}{}}

\begin{document}

\maketitle

\begin{abstract}
	Gagliardo--Nirenberg interpolation inequalities relate Lebesgue norms of iterated derivatives of a function.
    We present a generalization of these inequalities in which the low-order interpolation factor of the right-hand side is replaced by a Lebesgue norm of a pointwise product of derivatives of the function.
    We provide a short, independent proof, recovering the usual one-dimensional inequality as a particular case.
    We give examples of applications to control theory and PDE analysis, and discuss related open problems.
\end{abstract}

\section{Introduction}

The symbol $\lesssim$ denotes inequalities holding up to a constant, which can depend on the parameters of the statement, but not on the unknown function.
For convenience, we restrict the main statements to smooth functions, although they persist in appropriate Sobolev spaces, by standard regularization arguments.

\subsection{The classical Gagliardo--Nirenberg inequality}

In their full generality, the following interpolation inequalities date back to Gagliardo's and Nirenberg's respective short communications at the 1958 International Congress of Mathematicians in Edinburgh, later published in \cite{gag1,gag2,nir1,nir2}.
We refer to~\cite{fiorenza2021detailed} for a recent proof with a historical perspective.

\begin{theorem} \label{thm:GN}
	Let $p,q,r \in [1,\infty]$, $0 \leq k \leq j < m \in \N$ and $\theta \in [\theta^*,1]$ where
	\begin{equation}
		\theta^*:= \frac{j - k}{m - k}.
	\end{equation}
	Assume that
	\begin{equation}
		\label{eq:GN-relation}
		\frac{1}{p}-j = \theta \left( \frac{1}{r} - m \right)
		+ (1-\theta) \left(\frac{1}{q} - k\right).
	\end{equation}
	Then, for $u \in C^\infty_c(\R;\R)$,
	\begin{equation} \label{eq:GN-estimate}
		\| D^j u \|_{L^p(\R)} \lesssim \| D^m u \|_{L^r(\R)}^\theta \| D^k u \|_{L^q(\R)}^{1-\theta}.
	\end{equation}
\end{theorem}

In this 1D setting, estimate \eqref{eq:GN-estimate} had also been derived with optimal constants by Landau~\cite{landau}, Kolmogorov~\cite{kolmogoroff} and Stein~\cite{stein} for particular cases of the parameters.

\begin{remark}[Role of $k$]
	Usually, such inequalities are stated with $k = 0$.
	The less frequent case $k > 0$ can of course be reduced to the standard case $k = 0$ by applying the latter to the function $D^k u$.
	We include it in the above statement to highlight the symmetry with our generalization in \cref{thm:main}.
\end{remark}

\begin{remark}[Geometric setting] \label{rk:geometry}
	We restrict our statements to the 1D case.
	\emph{Mutatis mutandis} in \eqref{eq:GN-relation}, \cref{thm:GN} admits generalizations to~$\R^d$, or sufficiently regular bounded domains \cite{nir1}, or even exterior domains~\cite{crispo2004interpolation}, up to exceptional cases of the parameters (see e.g.\ \cite[Theorem 1.1 and below]{fiorenza2021detailed} for difficulties with $p = \infty$).
\end{remark}

\begin{remark}[Fractional versions] \label{rk:fractional}
	Although historical statements only involved integer orders of the derivatives, generalizations of \cref{thm:GN} to fractional Sobolev spaces~\cite{brezis2,brezis1} or Hölder spaces~\cite{kufner1995interpolation} are now well understood.
\end{remark}

\cref{thm:GN} extends to bounded domains, provided that one adds an (inhomogeneous) correction term, or a vanishing assumption, to rule out the case of low-degree polynomials (which could have $D^m u \equiv 0$ but $D^j u \neq 0$).

\begin{corollary}[Bounded domain] \label{cor:GN-bounded}
	Let $0 \leq k_0 \leq k$ and $s \in [1,\infty]$.
    
	For $u \in C^\infty([0,1];\R)$, without any vanishing assumption,
	\begin{equation} \label{eq:GN-bounded}
		\| D^j u \|_{L^p(0,1)} \lesssim \| D^m u \|_{L^r(0,1)}^\theta \| D^k u \|_{L^q(0,1)}^{1-\theta} + \| D^{k_0} u \|_{L^s(0,1)}.
	\end{equation}
    For $u \in C^\infty([0,1];\R)$ such that $u(0) = Du(0) = \dotsb = D^{m-1}u(0) = 0$,
	\begin{equation} \label{eq:GN-bounded-vanish}
		\| D^j u \|_{L^p(0,1)} \lesssim \| D^m u \|_{L^r(0,1)}^\theta \| D^k u \|_{L^q(0,1)}^{1-\theta}.
	\end{equation}
\end{corollary}

We will recover \cref{thm:GN} and \cref{cor:GN-bounded} as the particular case $\kappa = 1$ of our generalized versions \cref{thm:main} and \cref{cor:main-bounded}, of which we give self-contained proofs.

\subsection{Statement of the main result}

Motivated by applications to nonlinear control theory (see \cref{sec:control}) and PDE analysis (see \cref{sec:pde}), our main result is the following generalization of \cref{thm:GN} in which the low-order interpolation factor $\|D^k u\|_{L^q(\R)}$ of the right-hand side is replaced by a Lebesgue norm of a pointwise product of low-order derivatives of the function.

\begin{theorem} \label{thm:main}
	Let $p,q,r \in [1,\infty]$, $\kappa \in \N^*$ and 
 $0 \leq k_1 \leq \dotsb \leq k_{\kappa} \leq j < m \in \N$. 
	Let $\bar{k} := (k_1 + \dotsb + k_\kappa) / \kappa$.
	Let $\theta \in [\theta^*,1]$, where
	\begin{equation} \label{eq:theta*}
		\theta^* := \frac{j - \bar{k}}{m - \bar{k}}.
	\end{equation}
	Assume that
	\begin{equation} \label{eq:main-relation}
		\frac{1}{p}-j = \theta \left( \frac{1}{r} - m \right)
		+ (1-\theta) \left(\frac{1}{q \kappa} - \bar{k} \right).
	\end{equation}
	Then, for $u \in C^\infty_c(\R;\R)$,
	\begin{equation} \label{eq:main-estimate}
		\| D^j u \|_{L^p(\R)} \lesssim \| D^m u \|_{L^r(\R)}^\theta \| D^{k_1} u \dotsb D^{k_\kappa} u \|_{L^{q}(\R)}^{(1-\theta)/\kappa}.
	\end{equation}
\end{theorem}

\begin{remark} \label{rk:dkqk}
	Heuristically, everything behaves as if the pointwise product interpolation factor was replaced by $\| D^{\bar k} u \|_{L^{q \kappa}(\R)}$, see also \cref{sec:open}.
\end{remark}

\begin{remark} \label{rmk:pqr}
	In the critical case $\theta = \theta^*$, relation \eqref{eq:main-relation} is equivalent to
	\begin{equation} \label{eq:critic-pqr}
		\frac{1}{p} = \frac{\theta}{r} + \frac{1-\theta}{q \kappa}.
	\end{equation}
\end{remark}

\begin{corollary}[Bounded domain] 
    \label{cor:main-bounded}
	Let $0 \leq k_0 \leq k_1$ and $s \in [1,\infty]$.
    
	For $u \in C^\infty([0,1];\R)$, without any vanishing assumption,
    \begin{align}
        \label{eq:main-estimate-bounded}
		\| D^j u \|_{L^p(0,1)} & \lesssim \| D^m u \|_{L^r(0,1)}^\theta \| D^{k_1} u \dotsb D^{k_\kappa} u \|_{L^{q}(0,1)}^{(1-\theta)/\kappa} + \| D^{k_0} u \|_{L^s(0,1)}, \\
        \label{eq:main-estimate-bounded-bis}
		\| D^j u \|_{L^p(0,1)} & \lesssim \| D^m u \|_{L^r(0,1)}^\theta \| D^{k_1} u \dotsb D^{k_\kappa} u \|_{L^{q}(0,1)}^{(1-\theta)/\kappa} +\| D^{k_1} u \dotsb D^{k_\kappa} u \|_{L^{q}(0,1)}^{1/\kappa}.
    \end{align}
    For $u \in C^\infty([0,1];\R)$ such that $u(0) = Du(0) = \dotsb = D^{m-1}u(0) = 0$,
    \begin{equation}
        \label{eq:main-estimate-bounded-ter}
        \| D^j u \|_{L^p(0,1)} \lesssim \| D^m u \|_{L^r(0,1)}^\theta \| D^{k_1} u \dotsb D^{k_\kappa} u \|_{L^{q}(0,1)}^{(1-\theta)/\kappa}.
    \end{equation}
\end{corollary}

Our proof is inspired by Nirenberg's historical one, as rewritten recently in~\cite{fiorenza2021detailed}.
We combine an additive version of \eqref{eq:main-estimate} (see \cref{lem:compact}), a covering argument (see \cref{lem:covering}) and a summation argument (see \cref{sec:summation}).
Compared with the usual case, the main difference is in the proof of the additive version, where we tune the polynomial argument for the specific form of our low-order interpolation factor.

A related difficulty is that pointwise multiplicative inequalities of the form $|u'(x)|^2 \lesssim |u(x) u''(x)|$ usually require one to subtract from $u$ a local polynomial approximation, and to formulate the estimate using the Hardy--Littlewood maximal functions $Mu$, $Mu'$ and $Mu''$ instead of the raw functions (see e.g.\ \cite[Theorem 1]{js1994pointwise}).

We discuss related open problems in \cref{sec:open}.

\subsection{Some examples}

As illustrations of \cref{thm:main} for small values of the parameter $\kappa$, and in view of \cref{sec:control}, we state two particular cases.

\begin{corollary}
	Let $k \in \N^*$.
	For $u \in W^{2k,\infty}_0((0,1);\R)$,
	\begin{equation}
		\| D^k u \|_{L^6(0,1)}^6 \lesssim \| D^{2k} u \|_{L^\infty(0,1)}^2 \| u D^k u \|_{L^2(0,1)}^2.
	\end{equation}
\end{corollary}

\begin{proof}
	This follows from \cref{thm:main} with $\kappa = 2$, $\bar{k} = \frac k 2$, $\theta = \theta^* = \frac{k - \bar{k}}{2k - \bar{k}} = \frac 13$, for which $\frac 16 = \frac{\theta}{\infty} + \frac{1-\theta}{2 \cdot 2}$ so that \eqref{eq:critic-pqr} holds.
	The estimate for non-smooth $u$ follows by standard regularization arguments\footnote{\label{foot:Wm0}For $m \ge 1$, define $W^{m,\infty}_0((0,1);\R)$ as the space of functions $u \in W^{m,\infty}((0,1);\R)$ such that $u, Du, \dotsc, D^{m-1}u$ vanish at $0$ and $1$.
    Given such a $u$, its extension by $0$ belongs to $W^{m,\infty}(\R;\R)$.
    Then consider a family $(u_\varepsilon)_{\varepsilon \to 0}$ of mollified versions of this extension.
    For each $0 \le k < m$, $D^k u_\varepsilon \to D^k u$ strongly, while $\| D^m u_\varepsilon \|_{L^\infty} \leq \| D^m u \|_{L^\infty}$.
    Apply \cref{thm:main} to $u$ and pass to the limit.
    }.
\end{proof}

\begin{corollary} \label{cor:cubic}
	For $u \in W^{3,\infty}_0((0,1);\R)$,
	\begin{equation}
		\| u'' \|_{L^{12}(0,1)}^{12} \lesssim \| u''' \|_{L^\infty(0,1)}^6 \| u u' u'' \|_{L^2(0,1)}^2.
	\end{equation}
\end{corollary}

\begin{proof}
	This follows from \cref{thm:main} with $\kappa = 3$, $\bar{k} = 1$, $\theta = \theta^* = \frac{2-1}{3-1} = \frac 12$, for which $\frac 1 {12} = \frac \theta \infty + \frac{1-\theta}{2 \cdot 3}$ so that \eqref{eq:critic-pqr} holds. 
	The estimate for non-smooth $u$ follows by standard regularization arguments\footref{foot:Wm0}.
\end{proof}

\section{Proof of the main theorem}
\label{sec:proof}

\subsection{Additive Sobolev inequality with a product of derivatives}

We prove an additive version of \eqref{eq:main-estimate} on the interval $(0,1)$.
We denote by $\mathcal{P}_d$ the finite-dimensional space of polynomials of degree at most $d$, and by $\mu$ the Lebesgue measure on $\R$.
We start with the following sublevel-set version of the Remez inequality \cite{Remes1936}.

\begin{lemma} \label{lem:sublevel}
	Let $d \in \N$ and $\eta > 0$.
	There exists $c > 0$ such that, for all $Q \in \mathcal{P}_d$,
	\begin{equation} \label{eq:sublevel}
		\mu \left\{ x \in (0,1) : |Q(x)| < c \|Q \|_{L^\infty(0,1)} \right\} \leq \eta.
	\end{equation}
\end{lemma}

\begin{proof}
	By contradiction, for each $n \in \N^*$, there would exist $Q_n \in \mathcal{P}_d$ with $\|Q_n\|_{L^\infty} = 1$ for which~\eqref{eq:sublevel} fails with $c = 1/n$.
	Up to a subsequence, $Q_n \to Q$ uniformly on $[0,1]$, with $Q \neq 0$.
	For every $\delta > 0$ and $n$ large enough, $\{ |Q_n| < 1/n \} \subseteq \{ |Q| < \delta \}$ within $[0,1]$, so that $\mu \left\{ x \in (0,1) : |Q(x)| < \delta \right\} \geq \eta$.
	Letting $\delta \to 0$ yields $\mu \left\{ x \in (0,1) : Q(x) = 0 \right\} \geq \eta$, contradicting $Q \neq 0$.
\end{proof}

\begin{proposition} \label{lem:compact}
	Let $p,q,r \in [1,\infty]$, $\kappa \in \N^*$ and $0 \leq k_1 \leq \dotsb \leq k_{\kappa} \leq j < m \in \N$.
    
	For $u \in C^\infty([0,1];\R)$,
	\begin{equation} \label{eq:lem:compact}
		\| D^j u \|_{L^p(0,1)} 
		\lesssim
		\| D^m u \|_{L^r(0,1)} + \| D^{k_1}u \dotsb D^{k_\kappa} u \|_{L^q(0,1)}^{\frac 1\kappa}.
	\end{equation}
\end{proposition}

\begin{proof}
	By monotonicity of the Lebesgue spaces, we can take $p = \infty$ and $q = r = 1$.

    \medskip \noindent
    \emph{Decomposition of $u$.}
	Let $P \in \mathcal{P}_{m-1}$ be the Taylor polynomial of $u$ at $0$ of order $m-1$ and $R := u - P$.
	By Taylor's formula with integral remainder, for all $0 \leq s < m$,
	\begin{equation} \label{eq:taylor-remainder}
		\|D^s R\|_{L^\infty} \leq \| D^m u \|_{L^1},
	\end{equation}
	so that, by the triangle inequality,
	\begin{equation} \label{eq:Dju-DjP}
		\|D^j u\|_{L^\infty} \leq \|D^j P\|_{L^\infty} + \| D^m u \|_{L^1}.
	\end{equation}

    \medskip \noindent
	\emph{A good set.}
	Let $c > 0$ be given by \cref{lem:sublevel} with $d = m-1$ and $\eta = 1/(2\kappa)$, and let
	\begin{equation}
		F := \left\{ x \in (0,1) : |D^{k_i} P(x)| \geq c \|D^{k_i} P\|_{L^\infty} \text{ for all } 1 \leq i \leq \kappa \right\}.
	\end{equation}
	Applying \eqref{eq:sublevel} to each $D^{k_i} P \in \mathcal{P}_{m-1}$ and summing over $i$ gives $\mu(F) \geq 1/2$.
	
    On the finite-dimensional space $\mathcal{P}_{m-1}$, $\|\cdot\|_{L^\infty}$ is a norm and $D$ is continuous, so there exists $C \geq 1$, depending only on $m$, such that $\|DQ\|_{L^\infty} \leq C \|Q\|_{L^\infty}$ for all $Q \in \mathcal{P}_{m-1}$.
	Since $k_i \leq j$, writing $D^j P = D^{j-k_i}(D^{k_i} P)$ yields $\|D^j P\|_{L^\infty} \leq C^{j} \|D^{k_i} P\|_{L^\infty}$.
	Hence, setting $\rho := c C^{-j} > 0$, which depends only on $m$, $j$ and $\kappa$,
	\begin{equation} \label{eq:good-set}
		|D^{k_i} P(x)| \geq \rho \|D^j P\|_{L^\infty}
		\quad \text{for all } x \in F \text{ and } 1 \leq i \leq \kappa.
	\end{equation}

    \medskip \noindent
	\emph{Conclusion.}
	If $\rho \|D^j P\|_{L^\infty} \leq 2 \|D^m u\|_{L^1}$, then \eqref{eq:Dju-DjP} gives $\| D^j u \|_{L^\infty} \leq (1+\tfrac{2}{\rho}) \|D^m u\|_{L^1}$.
	Otherwise, for $x \in F$ and $1 \leq i \leq \kappa$, by \eqref{eq:good-set} and \eqref{eq:taylor-remainder} (legitimate since $k_i \leq j < m$),
	\begin{equation}
		|D^{k_i} u(x)| 
		\geq |D^{k_i} P(x)| - |D^{k_i} R(x)| 
		\geq \rho \| D^j P \|_{L^\infty} - \|D^m u\|_{L^1} 
		\geq \tfrac{\rho}{2} \|D^j P\|_{L^\infty}.
	\end{equation}
	Integrating over $F$ gives $\|D^{k_1}u \dotsb D^{k_\kappa}u\|_{L^1} \geq \tfrac 12 \left(\tfrac \rho2 \|D^j P\|_{L^\infty}\right)^{\kappa}$.
    Combining with \eqref{eq:Dju-DjP} yields $\|D^j u\|_{L^\infty} \leq \|D^m u\|_{L^1} + \tfrac{4}{\rho} \|D^{k_1}u \dotsb D^{k_\kappa} u\|_{L^1}^{1/\kappa}$.
\end{proof}

\begin{remark}
    \label{rk:ki<=j}
	The assumption $k_i \leq j$ is used above in the construction of the set~$F$, where it allows one to bound $\|D^j P\|_{L^\infty(0,1)}$ by $\|D^{k_i} P\|_{L^\infty(0,1)}$ for polynomials $P$ of degree at most $m-1$.
	This comparison fails when $k_i > j$, as does \eqref{eq:lem:compact} itself.
\end{remark}

\begin{corollary} \label{lem:scale}
	Let $p,q,r \in [1,\infty]$, $\kappa \in \N^*$ and $0 \leq k_1 \leq \dotsb \leq k_{\kappa} \leq j < m \in \N$.
    For any bounded interval $I \subset \R$ of positive length and $u \in C^\infty(\bar{I};\R)$,
	\begin{equation} \label{eq:lem:scale}
		|I|^{j - \frac 1 p} \| D^j u \|_{L^p(I)} 
		\lesssim
		|I|^{m - \frac 1 r} \| D^m u \|_{L^r(I)} 
		+ 
		|I|^{\bar{k} - \frac {1}{q\kappa}}
		\| D^{k_1}u \dotsb D^{k_\kappa} u \|_{L^q(I)}^{\frac 1\kappa},
	\end{equation}
	where $\bar{k} := (k_1 + \dotsb + k_\kappa) / \kappa$.
\end{corollary}

\begin{proof}
	This reduces to \cref{lem:compact} by a scaling argument.
	Write $I = (x_0, x_0 + |I|)$ for some $x_0 \in \R$.
	For $u \in C^\infty(\bar{I};\R)$, let $v \in C^\infty([0,1];\R)$ be defined by $v(x) := u(x_0 + x |I|)$, so that \eqref{eq:lem:scale} follows from \eqref{eq:lem:compact} with a constant independent of $|I|$.
\end{proof}

\subsection{Covering argument}
\label{sec:covering}

Let $u \in C^\infty([0,1];\R)$. 
We assume that $u$ is ``high-frequency'' in the sense of \eqref{eq:Dkv<Dmu} below, an assumption which is automatically satisfied when $u \in C^\infty_c((0,1);\R)$ (see \cref{sec:proof-main}).
We prove that we can find a covering of the support of $D^j u$ such that, on each interval, both terms of the right-hand side of \eqref{eq:lem:scale} are equal.
The proof is inspired by \cite[Lemma~3.3]{fiorenza2021detailed} and relies on Besicovitch's covering theorem \cite{besicovitch1945general}.

\begin{lemma}[Besicovitch] \label{lem:besicovitch}
	Let $E$ be a bounded subset of $\R$ and $r : E \to (0,+\infty)$.
	For $x \in E$, consider the open interval $I_x := (x-r_x,x+r_x)$.
	There exists a countable (possibly finite) collection of points $x_n \in E$ such that $E \subseteq \cup_n I_{x_n}$ and $\sum_n \mathbf{1}_{I_{x_n}} \leq 4$ on $\R$.
\end{lemma}

\begin{proof}
	This version corresponds to the 1D case of  \cite[Theorem~18.1, Chapter~2]{dibenedetto2002real}.
\end{proof}

\begin{proposition} \label{lem:covering}
	Let $q,r \in [1,\infty]$, $\kappa \in \N^*$ and $0 \leq k_1 \leq \dotsb \leq k_{\kappa} \leq j < m \in \N$.
    Let $\bar{k} := (k_1 + \dotsb + k_\kappa) / \kappa$.
	Let $u \in C^\infty([0,1];\R)$ and $v := D^{k_1} u \dotsb D^{k_\kappa} u$.
    Assume that
    \begin{equation} \label{eq:Dkv<Dmu}
        \| v \|_{L^q(0,1)}^{\frac 1 \kappa} 
        \leq 
        \| D^m u \|_{L^r(0,1)}.
    \end{equation}
    Define
	\begin{equation}
        \label{eq:def-E}
		E := \left\{ x \in (0,1) : D^j u(x) \neq 0 \right\}.
	\end{equation}
	There exists a countable family $(I_n)_n$ of non-empty open intervals $I_n \subset (0,1)$ such that
	\begin{align} \label{eq:lem:Ik-low}
		& E \subseteq \bigcup_n I_n, \\
		\label{eq:lem:Ik-high}
		& \sum_n \mathbf{1}_{I_n} \leq 4
		\quad
		\text{ on } \R,
	\end{align}
	and, for every $n$,
	\begin{equation} \label{eq:lem:Ik-2}
		|I_n|^{m - \frac 1 r} \| D^m u \|_{L^r(I_n)} 
		=
		|I_n|^{\bar{k} - \frac {1}{q\kappa}}
		\| D^{k_1}u \dotsb D^{k_\kappa} u \|_{L^q(I_n)}^{\frac 1\kappa}.
	\end{equation}
\end{proposition}

\begin{proof}
    For $x \in E$ and $h > 0$, set $J_x(h) := (x-h,x+h) \cap (0,1)$,
	\begin{align}
		\alpha_x(h) & := |J_x(h)|^{\bar{k} - \frac{1}{q\kappa}} \| v \|^{\frac 1\kappa}_{L^q(J_x(h))}, \\
		\beta_x(h) & := |J_x(h)|^{m-\frac{1}{r}} \|D^m u \|_{L^r(J_x(h))}.
	\end{align}
    On the one hand, for $h$ small enough, $J_x(h) = (x-h,x+h)$ and $|J_x(h)| = 2h$.
    In particular, as $h \to 0$, $\beta_x(h) = O(h^m)$.
    Moreover, by \cref{lem:scale} with $p = \infty$, $(2h)^j |D^j u(x)| \leq C (\alpha_x(h) + \beta_x(h))$.
    Since $x \in E$, $D^j u(x) \neq 0$.
    Since $j < m$, $\alpha_x(h) \gtrsim h^j$ and $\alpha_x(h) > \beta_x(h)$ for $h$ small enough. 
    
    On the other hand, for $h \geq \max \{ x, 1 - x \}$, $J_x(h) = (0,1)$ so that $\alpha_x(h) = \| v \|^{\frac 1\kappa}_{L^q(0,1)}$ and $\beta_x(h) = \|D^m u\|_{L^r(0,1)}$.
    By \eqref{eq:Dkv<Dmu}, $\alpha_x(h) \leq \beta_x(h)$.
    
    Hence, we can define
	\begin{equation} \label{eq:rx}
		r_x := \inf \{ h > 0 \mid \alpha_x(h) \leq \beta_x(h) \} \in (0,+\infty).
	\end{equation}
	In particular, by continuity, for every $x \in E$, $\alpha_x(r_x) = \beta_x(r_x)$.
    Consider the intervals $I_x := (x-r_x,x+r_x)$.
	By \cref{lem:besicovitch}, there exists a countable collection of points $x_n \in E$ such that $E \subseteq \cup_n I_{x_n}$ and $\sum_n \mathbf{1}_{I_{x_n}} \leq 4$ on $\R$.
    Let $I_n := I_{x_n} \cap (0,1) = J_{x_n}(r_{x_n})$.
    Then \eqref{eq:lem:Ik-low} and \eqref{eq:lem:Ik-high} are preserved and, on each $I_n$, \eqref{eq:lem:Ik-2} holds by definition of~$r_{x_n}$.
\end{proof}

\subsection{Proof of the high-frequency estimate}
\label{sec:summation}

We combine the additive Sobolev inequality and the covering argument to prove our main estimate in the high-frequency case.

\begin{proposition}
    \label{prop:high-freq}
    Under the assumptions of \cref{thm:main}, for $u \in C^\infty([0,1];\R)$ such that assumption \eqref{eq:Dkv<Dmu} holds, one has
    \begin{equation} 
        \label{eq:bound-high}
        \| D^j u \|_{L^p(0,1)} \lesssim \| D^m u \|_{L^r(0,1)}^\theta \| D^{k_1} u \dotsb D^{k_\kappa} u \|_{L^{q}(0,1)}^{(1-\theta)/\kappa}.
    \end{equation}
\end{proposition}

\begin{proof}
    Let $(I_n)_n$ be given by \cref{lem:covering} and set
    \begin{equation}
        a_n := \| D^m u \|_{L^r(I_n)}, \quad 
        b_n := \| D^{k_1}u \dotsb D^{k_\kappa}u \|_{L^q(I_n)}^{1/\kappa}, \quad 
        c_n := \| D^j u \|_{L^p(I_n)}.
    \end{equation}
    By the definition \eqref{eq:def-E} of~$E$ and \eqref{eq:lem:Ik-low},
    \begin{equation} 
        \label{eq:lp-sumislp}
        \| D^j u \|_{L^p(0,1)}
        = \| D^j u \|_{L^p(E)}
        \leq \| (c_n)_n \|_{\ell^p}.
    \end{equation}
    Using \cref{lem:scale}, there exists $C$ (independent of $u$) such that, for each $n$,
    \begin{equation}
        \label{eq:local-cn}
    	\begin{split}
    		c_n
    		& \leq C |I_n|^{\frac1p-j}
    		\big( |I_n|^{m - \frac 1 r} a_n
    		+ |I_n|^{\bar{k} - \frac {1}{q\kappa}} b_n\big) \\
    		& = 2 C |I_n|^{\frac1p-j} \big(|I_n|^{m - \frac 1 r} a_n \big)^{\theta}
    		\big(|I_n|^{\bar{k} - \frac {1}{q\kappa}}
    		b_n \big)^{1-\theta} 
            = 2C a_n^{\theta} b_n^{1-\theta}
    	\end{split}
    \end{equation}
    the second equality following from \eqref{eq:lem:Ik-2}, and the third from the scaling relation \eqref{eq:main-relation}.
    Define $p^* \in [1,\infty]$ by
    \begin{equation}
        \frac{1}{p^*} = \frac{\theta}{r} + \frac{1-\theta}{q \kappa}
        \quad \text{so that} \quad
        \frac{1}{p^*} - \frac{1}{p} = (m - \bar{k}) (\theta - \theta^*),
    \end{equation}
    using \eqref{eq:main-relation} and \eqref{eq:theta*}.
    Since $\theta \geq \theta^*$, $p^* \leq p$.
    Thus, using \eqref{eq:local-cn} and Hölder's inequality,
    \begin{equation}
        \label{eq:cn-holder-anbn}
        \| (c_n)_n \|_{\ell^p}
        \le \| (c_n)_n \|_{\ell^{p^*}}
        \le 2 C \| (a_n^\theta b_n^{1-\theta})_n \|_{\ell^{p^*}}
        \le 2 C \| (a_n)_n \|_{\ell^r}^\theta \| (b_n)_n \|_{\ell^{q\kappa}}^{1-\theta}.
    \end{equation}
    Finally, \eqref{eq:lem:Ik-high} gives
    \begin{equation}
        \label{eq:an-bn-norms}
        \|(a_n)_n\|_{\ell^r}\le4^{1/r}\|D^m u\|_{L^r(0,1)}
        \quad \text{and} \quad
        \|(b_n)_n\|_{\ell^{q\kappa}}
        \le4^{1/(q\kappa)}\|v\|_{L^q(0,1)}^{1/\kappa}.
    \end{equation}
    Gathering \eqref{eq:lp-sumislp}, \eqref{eq:cn-holder-anbn} and \eqref{eq:an-bn-norms} yields precisely \eqref{eq:bound-high}.
\end{proof}

\subsection{Proof of the theorem on the line}
\label{sec:proof-main}

We conclude the proof of \cref{thm:main}.
Since estimate \eqref{eq:main-estimate} is invariant under translations and rescalings, we may assume that $u \in C^\infty_c((0,1);\R)$.
Since all derivatives of $u$ vanish at $0$, by Taylor's formula with integral remainder, for each $1 \leq i \leq \kappa$, since $k_i < m$,
\begin{equation}
    \| D^{k_i} u \|_{L^\infty(0,1)} \leq \| D^m u \|_{L^1(0,1)},
\end{equation}
with unit constant since we are working on $(0,1)$.
By monotonicity of the Lebesgue norms on $(0,1)$, this entails that the high-frequency assumption \eqref{eq:Dkv<Dmu} holds.
Thus \cref{prop:high-freq} applies and concludes the proof.

\subsection{Proof of the corollary on the segment}

Under the assumptions of \cref{thm:main}, let $u \in C^\infty([0,1];\R)$ and $v := D^{k_1} u \dotsb D^{k_\kappa} u$.
We prove \cref{cor:main-bounded}.
If $u$ is high-frequency, in the sense that \eqref{eq:Dkv<Dmu} holds, \cref{prop:high-freq} applies and all the estimates of \cref{cor:main-bounded} hold (even without the additional correction terms).
We now assume conversely that $u$ is ``low-frequency'' in the sense that
\begin{equation}  
    \label{eq:low-freq}
    \| D^m u \|_{L^r(0,1)} < \| v \|_{L^q(0,1)}^{\frac 1\kappa}.
\end{equation}
The construction of the intervals of \cref{lem:covering} fails, but \cref{lem:compact} yields
\begin{equation} \label{eq:low-add}
    \| D^j u \|_{L^p(0,1)} \lesssim \| D^m u \|_{L^r(0,1)} + \| v \|_{L^q(0,1)}^{\frac 1\kappa}.
\end{equation}

\paragraph{Proof of \eqref{eq:main-estimate-bounded}.}
Since $k_0 \leq k_i < m$ for each $1 \leq i \leq \kappa$, \cref{lem:compact} with $\kappa = 1$ gives $\| D^{k_i} u \|_{L^\infty(0,1)} \lesssim \| D^m u \|_{L^r(0,1)} + \| D^{k_0} u \|_{L^s(0,1)}$.
Hence
\begin{equation} \label{eq:low-sob}
    \| v \|_{L^q(0,1)}^{\frac 1\kappa}
    \leq \prod_{i=1}^{\kappa} \| D^{k_i} u \|_{L^\infty(0,1)}^{\frac 1\kappa}
    \lesssim \| D^m u \|_{L^r(0,1)} + \| D^{k_0} u \|_{L^s(0,1)}.
\end{equation}
Moreover, by \eqref{eq:low-freq} and since $\theta \in [0,1]$,
\begin{equation}
    \label{eq:dmtheta}
    \| D^m u \|_{L^r(0,1)} \leq \| D^m u \|_{L^r(0,1)}^\theta \| v \|_{L^q(0,1)}^{\frac{1-\theta}{\kappa}}.
\end{equation}
Combining this with \eqref{eq:low-sob}, we obtain
\begin{equation}
    \| v \|_{L^q(0,1)}^{\frac 1\kappa}
    \lesssim
    \| D^m u \|_{L^r(0,1)}^\theta
    \| v \|_{L^q(0,1)}^{\frac{1-\theta}{\kappa}}
    + \| D^{k_0} u \|_{L^s(0,1)}.
\end{equation}
The conclusion follows from \eqref{eq:low-add} and \eqref{eq:dmtheta}.

\paragraph{Proof of \eqref{eq:main-estimate-bounded-bis}.}
Under \eqref{eq:low-freq}, estimate \eqref{eq:low-add} gives $\| D^j u \|_{L^p(0,1)} \lesssim \| v \|_{L^q(0,1)}^{1/\kappa}$.

\paragraph{Proof of \eqref{eq:main-estimate-bounded-ter}.}
When $u(0) = Du(0) = \dotsb = D^{m-1} u(0) = 0$, writing each $D^{k_i} u$ as an iterated integral of $D^m u$ gives $\| D^{k_i} u \|_{L^\infty(0,1)} \leq \| D^m u \|_{L^1(0,1)}$, whence
\begin{equation} \label{eq:low-vanish}
    \| v \|_{L^q(0,1)}^{\frac 1\kappa} \lesssim \| D^m u \|_{L^r(0,1)}.
\end{equation}
Combining \eqref{eq:main-estimate-bounded-bis} and \eqref{eq:low-vanish} proves \eqref{eq:main-estimate-bounded-ter}.

\section{An application to control theory}
\label{sec:control}

Our initial motivation concerns obstructions to small-time local controllability for nonlinear finite-dimensional scalar-input control-affine systems.
It is known that such obstructions are linked with interpolation inequalities (see \cite{BeauchardMarbach2026}).

As an example, given $p \in \N^*$, consider the following system on $\R^4$:
\begin{equation} \label{syst}
	\begin{cases}
		\dot{x}_1 = w \\
		\dot{x}_2 = x_1 \\
		\dot{x}_3 = x_2 \\
		\dot{x}_4 = x_1^2 x_2^2 x_3^2 - x_1^p
	\end{cases}
\end{equation}
with initial condition $x(0) = 0$ where $w \in L^\infty((0,T);\R)$ is the control to be chosen.
We are interested in the following local property.

\begin{definition}
	We say that system \eqref{syst} is \emph{small-time locally controllable} when, for every $T,\eta > 0$, there exists $\delta > 0$ such that, for every $x^* \in \R^4$ with $|x^*| \leq \delta$, there exists $w \in L^\infty((0,T);\R)$ such that $\|w\|_{L^\infty(0,T)} \leq \eta$ and the associated solution to \eqref{syst} with initial condition $x(0) = 0$ satisfies $x(T) = x^*$.
\end{definition}

\begin{proposition}
	System \eqref{syst} is small-time locally controllable if and only if 
    \begin{equation}
        p \in \{ 3, 5, 7, 8, 9, 10, 11 \}.   
    \end{equation}
\end{proposition}

\begin{proof}
	Let $T > 0$. 
	Consider $w \in L^\infty((0,T);\R)$ a control such that $x_1(T) = x_2(T) = x_3(T) = 0$.
    Then\footref{foot:Wm0} $u := x_3 \in W^{3,\infty}_0((0,T);\R)$ and 
	\begin{equation} \label{eq:x4}
		x_4(T) = \int_0^T (u u' u'')^2 - \int_0^T (u'')^p
	\end{equation}
	so that the possibility to reach a target of the form $(0,0,0,\pm 1)$ is linked with functional inequalities involving products of derivatives.
	We study each case.
	\begin{itemize}
		\item Case $p \geq 12$.
		First,
		\begin{equation}
			\| u'' \|_{L^p(0,T)}^p 
			\leq \| u'' \|_{L^\infty(0,T)}^{p-12} \| u'' \|_{L^{12}(0,T)}^{12}
		\end{equation}
		Moreover, since $u''' = x_3''' = w$ and $u''(0) = x_1(0) = 0$,
		\begin{equation}
			\| u'' \|_{L^\infty(0,T)} \leq T \| w \|_{L^\infty(0,T)}.
		\end{equation}
		Thus, thanks to the (scale-invariant) interpolation inequality of \cref{cor:cubic},
		\begin{equation}
			\| u'' \|_{L^p(0,T)}^p 
			\leq C T^{p-12} \| w \|_{L^\infty(0,T)}^{p-6} \int_0^T (u u' u'')^2.
		\end{equation}
		Substituting in \eqref{eq:x4} proves that $x_4(T) \geq 0$ when $C T^{p-12} \| w \|_{L^\infty(0,T)}^{p-6} \leq 1$.
		Thus, choosing $0 < \eta \ll 1$ such that $C T^{p-12} \eta^{p-6} \leq 1$ negates the definition of small-time local controllability.
		
		\item Case $7 \leq p \leq 11$.
		Let $0 \neq \chi \in C^\infty_c((0,T);\R)$ and consider $w(t) := \varepsilon \chi'''(t)$ for $0 < \varepsilon \ll 1$.
		As $\varepsilon \to 0$, $w \to 0$ in $L^\infty((0,T);\R)$.
		Moreover, by \eqref{eq:x4},
		\begin{equation}
			x_4(T) = \varepsilon^6 \int_0^T (\chi \chi' \chi'')^2 + O(\varepsilon^7).
		\end{equation}
		So one can move in the direction $(0,0,0,+1)$.
		
		Conversely, set $u(t) = \varepsilon^{1+3a} \chi(t\varepsilon^{-a})$ or equivalently $w(t) := \varepsilon \chi'''(t \varepsilon^{-a})$ for $a > 0$ and $0 < \varepsilon \ll 1$.
		As $\varepsilon \to 0$, $w \to 0$ in $L^\infty((0,T);\R)$.
		By \eqref{eq:x4},
		\begin{equation}
			x_4(T) = \varepsilon^{6+13a} \int_0^T (\chi \chi' \chi'')^2 - \varepsilon^{p+(p+1)a} \int_0^T (\chi'')^p.
		\end{equation}
		If $\int_0^T (\chi'')^p > 0$ and $p(1+a) < 6+12a$ (which is possible when $7 \leq p \leq 11$), one can move in the direction $(0,0,0,-1)$.
		
		From these elementary movements, it is classical to conclude that \eqref{syst} is small-time locally controllable.
		
		\item Case $p = 1$. 
		Then $\dot{x}_2 + \dot{x}_4 = (x_1 x_2 x_3)^2 \geq 0$.
		Hence, for any control, $(x_2 + x_4)(T) \geq 0$ so targets with $x_2^* + x_4^* < 0$ are not reachable.
		
		\item Case $p \in \{ 2, 4, 6 \}$.
		The system does not satisfy Stefani's necessary condition for small-time local controllability (see \cite{stefani}).
		
		\item Case $p \in \{ 3, 5 \}$.
		The system satisfies Hermes' sufficient condition for small-time local controllability (see \cite{hermes1982}).\qedhere
	\end{itemize}
\end{proof}

\section{An application to PDE analysis}
\label{sec:pde}

Our inequalities are also relevant to the analysis of nonlinear PDEs, where pointwise products of derivatives arise as dissipation terms (see e.g.\ \cite{BerettaBertschDalPasso1995,Bernis1996} for thin-film models).
\cref{thm:main} allows one to retain them as low-order interpolation factors.
We illustrate this on a one-dimensional regularization of the porous medium equation.

Let $I := (0,1)$, $T > 0$, $Q_T := (0,T) \times I$ and $\varepsilon \in (0,1]$.
We consider
\begin{equation} \label{eq:CH-system}
    \begin{cases}
        \partial_t u + \varepsilon \partial_x^4 u - \partial_x^2 (u^3) = 0 & \text{in } Q_T, \\
        \partial_x u = \partial_x^3 u = 0 & \text{on } (0,T) \times \partial I, \\
        u(0) = u_0 & \text{in } I,
    \end{cases}
\end{equation}
which converges formally, as $\varepsilon \to 0$, to the porous medium equation $\partial_t u = \partial_x^2 (u^3)$.

Such limits are investigated, in higher dimension, by Colli and Fukao in~\cite{ColliFukao2016}, where nonlinear diffusion equations are obtained from Cahn--Hilliard systems; \eqref{eq:CH-system} is a one-dimensional counterpart of their porous medium example~\cite[Example~2]{ColliFukao2016}, with $q = 3$.
 
Since our purpose is to derive \emph{a priori} estimates, we only consider smooth solutions, that is, functions $u \in C^\infty([0,T] \times \bar{I};\R)$ satisfying \eqref{eq:CH-system} pointwise.
Throughout this section, $\lesssim$ denotes an inequality up to a constant depending on $T$ and $u_0$, but neither on $\varepsilon$ nor on $u$.

\begin{lemma} \label{lem:CH-energy}
    Smooth solutions of \eqref{eq:CH-system} satisfy
    \begin{equation} \label{eq:CH-bounds}
        \|u\|_{L^\infty((0,T);L^2(I))}
        + \sqrt{\varepsilon} \, \|\partial_x^2 u\|_{L^2(Q_T)}
        + \|u \partial_x u\|_{L^2(Q_T)} \lesssim 1.
    \end{equation}
\end{lemma}

\begin{proof}
    We multiply \eqref{eq:CH-system} by $u$ and integrate over $I$.
    Since $\partial_x u$ and $\partial_x^3 u$ vanish on $\partial I$, so does $\partial_x(u^3) = 3 u^2 \partial_x u$, and all boundary terms vanish.
    Hence $\int_I u \partial_x^4 u = \|\partial_x^2 u\|_{L^2(I)}^2$ and $- \int_I u \partial_x^2 (u^3) = 3 \|u \partial_x u\|_{L^2(I)}^2$, which gives the energy identity
    \begin{equation}
        \frac 12 \frac{\dd}{\dd t} \|u\|_{L^2(I)}^2
        + \varepsilon \|\partial_x^2 u\|_{L^2(I)}^2
        + 3 \|u \partial_x u\|_{L^2(I)}^2 = 0
    \end{equation}
    Integrating it in time yields \eqref{eq:CH-bounds}.
\end{proof}

The nonlinear dissipation is a pointwise product of derivatives of orders $0$ and $1$.
It can therefore be used directly as the low-order interpolation factor of \cref{cor:main-bounded}.

\begin{proposition} \label{prop:CH-product}
    Uniformly for $\varepsilon \in (0,1]$, smooth solutions of \eqref{eq:CH-system} satisfy
    \begin{equation} \label{eq:CH-product}
        \varepsilon^{1/6} \|\partial_x u\|_{L^3(Q_T)} \lesssim 1.
    \end{equation}
\end{proposition}

\begin{proof}
    We apply \eqref{eq:main-estimate-bounded} with $\kappa = 2$, $(k_1,k_2) = (0,1)$, $j = 1$, $m = 2$, $q = r = 2$, $k_0 = 0$ and $s = 2$, for which $\bar{k} = \frac 12$, $\theta = \theta^* = \frac 13$ and $p = 3$.
    For every $t \in (0,T)$, we obtain\footnote{For the particular exponents involved here, the interpolation inequality $\| u' \|_{L^3}^3 \leq \| u'' \|_{L^2} \| u u' \|_{L^2}$ could also be obtained using integration by parts.}
    \begin{equation}
        \|\partial_x u\|_{L^3(I)}
        \lesssim \|\partial_x^2 u\|_{L^2(I)}^{1/3} \|u \partial_x u\|_{L^2(I)}^{1/3}
        + \|u\|_{L^2(I)}.
    \end{equation}
    Cubing, integrating in time and using the Cauchy--Schwarz inequality, then \eqref{eq:CH-bounds},
    \begin{equation}
        \|\partial_x u\|_{L^3(Q_T)}^3
        \lesssim \|\partial_x^2 u\|_{L^2(Q_T)} \|u \partial_x u\|_{L^2(Q_T)} + 1
        \lesssim \varepsilon^{-1/2},
    \end{equation}
    which is \eqref{eq:CH-product}.
\end{proof}

We now compare \eqref{eq:CH-product} with the estimate obtained from the same energy identity through the classical Gagliardo--Nirenberg inequality.

\begin{lemma} \label{lem:CH-linfty}
    Smooth solutions of \eqref{eq:CH-system} satisfy $\|u\|_{L^6((0,T);L^\infty(I))} \lesssim 1$.
\end{lemma}

\begin{proof}
    Let $v := u^2$, so that $\partial_x v = 2 u \partial_x u$ and $\|v\|_{L^1(I)} = \|u\|_{L^2(I)}^2$.
    By \eqref{eq:GN-bounded} applied to $v$, with $j = k = k_0 = 0$, $m = 1$, $p = \infty$, $q = s = 1$, $r = 2$ and $\theta = \frac 23$, raised to the power $3$,
    \begin{equation*}
        \|u\|_{L^\infty(I)}^6 = \|v\|_{L^\infty(I)}^3
        \lesssim \|u \partial_x u\|_{L^2(I)}^2 \|u\|_{L^2(I)}^2 + \|u\|_{L^2(I)}^6.
    \end{equation*}
    Integrating in time and using \eqref{eq:CH-bounds} concludes.
\end{proof}

\begin{proposition} \label{prop:CH-GN}
    Uniformly for $\varepsilon \in (0,1]$, smooth solutions of \eqref{eq:CH-system} satisfy
    \begin{equation} \label{eq:CH-GN}
        \varepsilon^{1/4} \|\partial_x u\|_{L^3(Q_T)} \lesssim 1.
    \end{equation}
\end{proposition}

\begin{proof}
    By \eqref{eq:GN-bounded} with $j = 1$, $m = 2$, $k = k_0 = 0$, $p = 3$, $q = 6$, $r = s = 2$ and $\theta = \theta^* = \frac 12$, for every $t \in (0,T)$,
    \begin{equation*}
        \|\partial_x u\|_{L^3(I)}
        \lesssim \|\partial_x^2 u\|_{L^2(I)}^{1/2} \|u\|_{L^6(I)}^{1/2} + \|u\|_{L^2(I)}.
    \end{equation*}
    Since $\frac 13 = \frac 12 \cdot \frac 12 + \frac 12 \cdot \frac 16$, Hölder's inequality in time gives
    \begin{equation*}
        \|\partial_x u\|_{L^3(Q_T)}
        \lesssim \|\partial_x^2 u\|_{L^2(Q_T)}^{1/2} \|u\|_{L^6(Q_T)}^{1/2}
        + T^{1/3} \|u\|_{L^\infty((0,T);L^2(I))},
    \end{equation*}
    and one concludes with \eqref{eq:CH-bounds} and \cref{lem:CH-linfty}.
\end{proof}

Hence, as $\varepsilon \to 0$, the estimate of \cref{prop:CH-product} obtained using our inequality is less degenerate than the one obtained in \cref{prop:CH-GN} using the usual inequality.

\begin{remark}
    In the regularized system of \cite{ColliFukao2016}, the analogue of the nonlinear dissipation is $\int_\Omega \beta_\lambda'(u_{\varepsilon,\lambda}) |\nabla u_{\varepsilon,\lambda}|^2$, which appears in \cite[Eq.~(6.2)]{ColliFukao2016} but is discarded, by monotonicity of $\beta_\lambda$, in the passage to \cite[Eq.~(6.5)]{ColliFukao2016}, and hence does not occur in \cite[Eq.~(6.6)]{ColliFukao2016}.
    For $\beta(r) = r^3$, it is exactly $3\|u \partial_x u\|_{L^2(I)}^2$.
    The bounds retained there, namely $u^3$ bounded in $L^2(Q_T)$ by \cite[Eq.~(2.18)]{ColliFukao2016} and $\sqrt{\varepsilon} \, \partial_x^2 u$ bounded in $L^2(Q_T)$ by \cite[Eq.~(6.6)]{ColliFukao2016}, are precisely the ones used in the proof of \cref{prop:CH-GN}, which only requires $\|u\|_{L^6(Q_T)} \lesssim 1$.
\end{remark}

\section{Related open problems}
\label{sec:open}

\subsection{Higher dimensions}

As noted in \cref{rk:geometry}, standard Gagliardo--Nirenberg inequalities admit generalizations in $\R^d$.
Extending \cref{thm:main} accordingly is a natural goal, though identifying the appropriate multidimensional analogue of the product $D^{k_1} u \dotsb D^{k_\kappa} u$ (in particular, which collections of partial derivatives or multiindices should enter), may be delicate.

\subsection{Fractional Sobolev spaces}

As noted in \cref{rk:fractional}, standard Gagliardo--Nirenberg inequalities extend to fractional Sobolev spaces. 
Naturally, one could seek similar generalizations of \cref{thm:main}.
As seen in \cref{rk:dkqk}, even for integer parameters, $D^{k_1} u \dotsb D^{k_\kappa} u$ already behaves as a fractional Sobolev norm when $\bar{k} \notin \N$.

\subsection{Two necessary assumptions}

In \cref{thm:main}, the relation \eqref{eq:main-relation} between the parameters corresponds to the scaling invariance by spatial dilations (replace $u(x)$ by $u(\lambda x)$).
It is therefore necessary.

Similarly the condition $\theta \geq \theta^*$ is forced by high-frequency oscillations.
Indeed, given $\varphi \in C^\infty_c(\R;\R)$, consider the family $\varphi(x) \cos (Nx)$ with $N \gg 1$.
The left-hand side is of order $N^j$ whereas the right-hand side is bounded above by $N^{\theta m} N^{(1-\theta) \bar{k}}$.
Thus, we must have $j \leq \theta m + (1-\theta) \bar{k}$, which, by \eqref{eq:theta*}, is equivalent to $\theta \ge \theta^*$.

\subsection{About the assumption \texorpdfstring{$k_i \leq j$}{ki <= j}}

While our proofs in \cref{sec:proof} use the condition $k_i \leq j$ for all $1 \leq i \leq \kappa$ from \cref{thm:main} (see \cref{rk:ki<=j}), examples below show that this assumption is not always necessary.

As we noted in \cref{rk:dkqk}, $\| D^{k_1} u \dotsb D^{k_\kappa} u \|_{L^q}^{1/\kappa}$ heuristically behaves as $\| D^{\bar{k}} u \|_{L^{q\kappa}}$.
Since we want this to be a low-order interpolation factor for $\| D^j u \|_{L^p}$, a natural weaker assumption is $\bar{k} \leq j$.
In fact, this is a consequence\footnote{Up to the exceptional end-point case ($\theta = 0$, $p = \infty$, $q = 1$, $\bar{k} = j + \tfrac 1 \kappa$).} of \eqref{eq:main-relation} when $k_i < m$ and $\theta \in [0,1]$.

To explore the validity of the assumption $\bar{k} \leq j$, one may more specifically ask when the following result holds (formally corresponding to $j = \bar{k}$ and $\theta = \theta^* = 0$).

\begin{open} \label{open}
	Let $q \in [1,\infty]$, $\kappa \in \N^*$ and $0 \leq k_1 \leq \dotsb \leq k_{\kappa} \in \N$. 
	Let $\bar{k} := (k_1 + \dotsb + k_\kappa) / \kappa$.
	When is it true that, for $u \in C^\infty_c(\R;\R)$,
	\begin{equation} \label{eq:dbarku}
		\| D^{\bar{k}} u \|_{L^{q \kappa}(\R)} \lesssim \| D^{k_1} u \dotsb D^{k_\kappa} u \|_{L^{q}(\R)}^{\frac 1 \kappa},
	\end{equation}
	where, when $\bar{k}$ is not an integer, the left-hand side should be interpreted as an appropriate fractional semi-norm of $u$.
\end{open}

Following the heuristic of \cref{rk:dkqk}, positive answers to \cref{open} could give the desired extension of \cref{thm:main} under the weaker assumption $\bar{k} \leq j$ (up to exceptional cases) thanks to the fractional Gagliardo--Nirenberg inequalities of the form $\| D^j u \|_{L^p(\R)} \lesssim \| D^m u \|_{L^r(\R)}^\theta \| D^{\bar k} u \|_{L^{q\kappa}(\R)}^{1-\theta}$ (see for instance \cite{brezis2,brezis1} when $\bar{k} \notin \N$).

As an example of this mechanism, the interpolation inequality of \cref{cor:cubic} can be alternatively derived from  the usual Gagliardo--Nirenberg estimate of \cref{thm:GN}
\begin{equation}
    \| u'' \|_{L^{12}(\R)} \lesssim \| u''' \|_{L^\infty(\R)}^{\frac 12} \| u' \|_{L^6(\R)}^{\frac 12}
\end{equation}
and the coercivity estimate \eqref{eq:uu'u''-mino} proved below.

\subsubsection{Some positive answers}

We give \emph{ad hoc} arguments showing that \eqref{eq:dbarku} holds in specific examples, suggesting that \cref{open} may hold for some classes of parameters.

\begin{lemma}
	For $u \in C^\infty_c(\R;\R)$,
	\begin{align}
		\label{eq:w12-4}
		\| u \|_{\dot{W}^{\frac 12, 4}(\R)} & \lesssim \| u u' \|^{\frac 12}_{L^2(\R)}, \\
		\label{eq:u'l4}
		\| u' \|_{L^4(\R)} & \lesssim \| u u'' \|_{L^2(\R)}^{\frac 12}, \\
		\label{eq:uu'u''-mino}
		\| u' \|_{L^6(\R)} & \lesssim \| u u' u'' \|_{L^2(\R)}^{\frac 13}.
	\end{align}
\end{lemma}

\begin{proof}
	Estimate \eqref{eq:w12-4} can be derived from the remark that $(u^2)' = 2 u u'$.
	Hence
	\begin{equation}
		\| u \|_{\dot{W}^{\frac 12, 4}(\R)}
		\lesssim \| |u| \|_{\dot{W}^{\frac 12, 4}(\R)}
		\lesssim \| u^2 \|_{\dot{H}^1(\R)}^{\frac 12}
		= \| 2 u u' \|_{L^2(\R)}^{\frac 12},
	\end{equation}
	where the first estimate with the absolute value is derived in \cite[Théorème 2]{lemarie} and the second estimate in \cite[Section 5.4.4]{runst1996sobolev}.

    Integrations by parts and the Cauchy--Schwarz inequality yield:
	\begin{align}
		\int_\R (u')^4 & = - 3 \int_\R u (u')^2 u''
		\leq 3 \left( \int_\R (u')^4 \right)^{\frac 12} \left( \int_\R (u u'')^2 \right)^{\frac 12}, \\
		\int_\R (u')^6 & = - 5 \int_\R u (u')^4 u''
		\leq 5 \left( \int_\R (u')^6 \right)^{\frac 12} \left( \int_\R (u u' u'')^2 \right)^{\frac 12},
	\end{align}
	which entail \eqref{eq:u'l4} and \eqref{eq:uu'u''-mino}.
\end{proof}

Estimate \eqref{eq:u'l4} above is very classical, for example stated as Lemma~1 in~\cite{kalamajska2012some}, which contains many interesting generalizations.

\subsubsection{Some negative answers}

Conversely, \eqref{eq:dbarku} fails for some parameters, so that a positive answer to \cref{open} must exclude at least the following three situations.

\begin{lemma} \label{lem:counter}
	None of the following estimates holds for all $u \in C^\infty_c(\R;\R)$:
	\begin{align}
		\label{eq:no-h12}
		\| u \|_{\dot{H}^{\frac 12}(\R)} & \lesssim \| u u' \|_{L^1(\R)}^{\frac 12}, \\
		\label{eq:no-u2u3}
		\| u' \|_{L^6(\R)} & \lesssim \| u^2 u''' \|_{L^2(\R)}^{\frac 13}, \\
		\label{eq:no-uu2}
		\| u' \|_{L^\infty(\R)} & \lesssim \| u u'' \|_{L^\infty(\R)}^{\frac 12}.
	\end{align}
\end{lemma}

\begin{proof}
    Fix $\eta \in C^\infty(\R;[0,1])$, nondecreasing, with $\eta = 0$ on $(-\infty,0]$ and $\eta = 1$ on $[1,+\infty)$.

    \medskip\noindent
	\emph{Proof for \eqref{eq:no-h12}.}
	For $0 < \varepsilon < \frac 14$, let $u(x) := \eta(x/\varepsilon) \, \eta((1-x)/\varepsilon)$.
	Since $u^2$ increases from $0$ to $1$, then decreases to $0$, one has $\| u u' \|_{L^1} = 1$.
	Moreover, $u = 1$ on $[\varepsilon,1-\varepsilon]$ and $u = 0$ on $(-\infty,0]$, so that, up to a constant, the Slobodeckij semi-norm satisfies
	\begin{equation}
		\| u \|_{\dot{H}^{\frac 12}}^2 
		= \iint_{\R^2} \frac{|u(x) - u(y)|^2}{|x-y|^2}
		\geq \int_{\varepsilon}^{1-\varepsilon} \int_{-\infty}^{0} \frac{\dd y \dd x}{(x-y)^2}
		= \log \frac{1-\varepsilon}{\varepsilon} \xrightarrow[\varepsilon \to 0]{} + \infty.
	\end{equation}
	
	\medskip \noindent
	\emph{Proof for \eqref{eq:no-u2u3}.}
	For $0 < \varepsilon < \frac 14$, let $u(x) := x(1-x) \eta(x/\varepsilon) \, \eta((1-x)/\varepsilon)$.
	On $[\varepsilon,1-\varepsilon]$, $u(x) = x(1-x)$, so $u''' = 0$ there.
	Hence $u^2 u'''$ is supported in two intervals of length $\varepsilon$, on which Leibniz' rule gives $|u| \lesssim \varepsilon$ and $|u'''| \lesssim \varepsilon^{-2}$.
	Thus $\| u^2 u''' \|_{L^2} \lesssim \varepsilon^{1/2} \to 0$, whereas $u'(x) = 1 - 2x$ on $[\frac 14, \frac 13]$, so that $\| u' \|_{L^6} \gtrsim 1$.
	
	\medskip \noindent
	\emph{Proof for \eqref{eq:no-uu2}.}
	Let $\phi \in C^\infty_c((-2,2);\R)$ with $\phi = 1$ on $[-1,1]$.
	For $n \geq 1$, we define $u(x) := x \phi\left(\frac{\log x}{n}\right)$ for $x > 0$ and $u(x) := 0$ for $x \leq 0$.
	Then $\supp u \subset [e^{-2n}, e^{2n}]$, so $u$ vanishes near $0$ and $u \in C^\infty_c(\R;\R)$.
	Writing $t := \frac{\log x}{n}$, one has, for $x > 0$,
	\begin{equation}
		u'(x) = \phi(t) + \tfrac 1n \phi'(t), 
		\qquad
		u''(x) = \tfrac 1{nx} \left( \phi'(t) + \tfrac 1n \phi''(t) \right),
	\end{equation}
	so that $u u'' = \frac 1n \phi(t) \left( \phi'(t) + \frac 1n \phi''(t) \right)$ and $\| u u'' \|_{L^\infty} \lesssim \frac 1n \to 0$.
    
	But $u(x) = x$ on $[e^{-n}, e^{n}]$, so $\| u' \|_{L^\infty} \geq 1$.
\end{proof}

\section*{Acknowledgements}

The author is deeply indebted to Karine Beauchard for numerous discussions and preliminary results on particular cases of \cref{thm:main}, as well as for the underlying control-theoretic motivation.
The author also thanks Frédéric Bernicot and Cristina Benea for encouraging discussions concerning \cref{open}.

The author is supported by ANR-24-CE40-4571.

\bibliographystyle{plain}
\bibliography{main}

\end{document}